\documentclass[journal]{IEEEtran}
\usepackage{cite}
\usepackage{caption}
\usepackage{makecell}
\usepackage{multirow}
\usepackage{mathtools}
\makeatletter
\newcommand*{\rom}[1]{\expandafter\@slowromancap\romannumeral #1@}
\makeatother
\usepackage[font=small]{caption}

\ifCLASSINFOpdf
   \usepackage[pdftex]{graphicx}
\else
   \usepackage[dvips]{graphicx}
\fi
\usepackage{graphicx}
\usepackage{amsmath}
\usepackage{algorithmic}
\usepackage{array}
\usepackage{placeins}
\ifCLASSOPTIONcompsoc
  \usepackage[caption=false,font=normalsize,labelfont=sf,textfont=sf]{subfig}
\else
  \usepackage[caption=false,font=footnotesize]{subfig}
\fi
\usepackage{url}
\begin{document}
\title{Integrating Ultra-Fast Charging Stations within the Power Grids of Smart Cities: A Review}

\author{Danielle~Meyer,~\IEEEmembership{Student Member,~IEEE,}
        and~Jiankang Wang,~\IEEEmembership{Member,~IEEE} \thanks{D. Meyer and JK Wang are with the Department of Electrical and Computer Engineering, The Ohio State University, Columbus, OH, 43235 USA (email: meyer.758@buckeyemail.osu.edu, wang.6536@osu.edu).}}%
        
{

\maketitle

\begin{abstract}
Plug-In Electric Vehicles (PEV) have become a key factor driving towards smart cities, which allow for higher energy efficiency and lower environmental impact across urban sectors. Industry vision for future PEV includes the ability to recharge a vehicle at a speed comparable to traditional gas refueling, i.e., less than $3$ minutes per vehicle. Such a technology, referred to as Ultra-Fast Charging (UFC), has drawn much interest from research and industry. The large power density, impulsive, and stochastic loading characteristics of UFC, however, pose unprecedented challenges to existing electricity supply infrastructure. Planning the locations and electric capacities of these UFC stations is critical to preventing detrimental impacts, including grid asset depreciation, grid instabilities, and deteriorated power quality. In this paper, we first review planning methods for conventional charging stations and then discuss outlooks for UFC planning solutions by drawing an analogy with renewable energy source planning, which presents similar power density and stochastic characteristics as UFC. While this paper mainly focuses on UFC planning from the power grid perspective, other urban aspects, including traffic flow and end-user behavior, are examined for feasible UFC integration within smart cities. 
\end{abstract}
\begin{IEEEkeywords}
Ultra-Fast Charging (UFC), Smart Cities, Renewable Integration, Traffic Flow, User Behavior, Charging Station (CS) Planning
\end{IEEEkeywords}

%
\IEEEpeerreviewmaketitle

\section{Introduction}
\IEEEPARstart{T}{he} past seven years have seen consistent growth in Plug-in Electric Vehicle(PEV) market share. In particular, the United States saw a 22\% increase in EV sales from 2015 to 2016 \cite{IEA_Outlook}. Country specific targets, shown in Table \ref{targets}, further illustrate the continued push towards increased PEV adoption. The motivation for these aggressive targets includes the numerous benefits of PEVs, which lend themselves to the vision of smart cities. 

One of these benefits is urban environmental improvement. According to the U.S. Energy Information Agency, $\text{CO}_2$ emissions from the transport sector have been higher than those of the other sectors, but these numbers begin to show a decline in areas with higher PEV penetration rates \cite{EIA_2017}. Further, integrating charging technology with renewable energy and storage can synergize with electricity supply infrastructure, i.e., the power grid, to increase environmental benefits and improve energy efficiency  \cite{EIA_2017}.  From autonomous vehicles to vehicle-to-grid (V2G), advanced PEV features can provide a sustainable use of resources and a higher quality of life. 

Charging stations (CS) have a critical role in facilitating increased PEV pentration. Residential charging, AC Level 1, can take up to $13$ hours for a full charge, while AC Levels 2 and 3 and DC charging can significantly reduce charging to times to within an hour, as shown in Table \ref{levels}. The latter methods are only available for commercial users and at public CS based on SAE standards \cite{sae2011sae}. Ultra-Fast Charging (UFC) has been researched by academia and industry, aiming to reach charging speeds comparable to traditional gas refueling. 

\begin{table*}[t!]
\begin{center}
\caption {Global Leading EV Market Status and Ambitions \cite{IEA_Outlook}} \label{targets} 
    \begin{tabular}{ | c | c | c | c | c | c |}
    \hline
    Country & \thead{New EV Registrations \\ in 2016 (thousands)} &\thead{Cumulative EV Sales \\ through 2016 (thousands)} & EV Market Share & Penetration Targets \\ \hline
    United States & 159.62 & 563.71 & 0.9\% & 10\% in 2025 \\ \hline
    China & 336.00 & 648.77 & 1.4\% & 20\% by 2025 \\ \hline
    Norway & 50.18 & 133.26 & 28.8\% & 100\% by 2025 \\ \hline
    Netherlands & 24.48 & 112.01 & 6.4\% & 100\% by 2035 \\ \hline
    United Kingdom & 37.91 & 86.42 & 1.4\% & 100\% by 2040 \\ \hline
    Japan & 24.85 & 151.25 & 0.6\% & 20\%$\sim$30\% by 2030 \\ \hline
    \multicolumn{5}{c}{Note 1: EV includes Plug-in Electric Vehicle (PEV) and Plug-in Hybrid Electric Vehicle (PHEV).}\\
    \multicolumn{5}{c}{Note 2: Battery is the sole power source for BEV and an engine is paired with a battery to power a PEV.}\\
    \multicolumn{5}{c}{Note 3: The market share is calculated on the basis of total market size of each country.}\\
    \end{tabular}
\end{center}
\end{table*}

\begin{table*}[!t]
\begin{center}
\caption{PEV Charging Methods, Configurations and Ratings \cite{IEA_Outlook,shareef2016review,sae2011sae}} \label{levels} 
    \begin{tabular}{ | c | c | c | c | c | c |}
    \hline
    Charging Method & \thead{Port Terminal Voltage \\ and Connection} & \thead{Maximum Power \\ (kW)} & Charging Time & System Level Connection \\ \hline
    AC Level 1 & 120V 1-phase & 2 & 10$\sim$13h & \multirow{3}{*}{\thead{Residential/Commercial \\ Secondary Customer \\ 120V and 240V}} \\ \cline{1-4}
    AC Level 2 & 240V 1-phase & 20 & 1$\sim$4h  &\\ \cline{1-4}
    AC Level 3 & 240V 3-phase & 43.5 & $\sim$1h & \\ \hline
    DC Level 1 & 200$\sim$450V & 36 & 0.5$\sim$1.44h & \multirow{3}{*}{\thead{Commercial \\ Primary Customer \\ 480V}} \\ \cline{1-4}
    DC Level 2 & 200$\sim$450V & 96 & 0.2$\sim$0.58h &  \\ \cline{1-4}
    DC Level 3 & 200$\sim$600V & 200 & $\sim$10min & \\ \hline
    Ultra-Fast Charging & $\geq$800V & 1,000 & $\sim$Gas Refueling & \thead{Sub-transmission \\ Primary Customer \\ 26kV or 69kV} \\ \hline
    \multicolumn{5}{c}{Note1: Ultra-Fast Charging is not yet finalized }\\
    \end{tabular}
\end{center}
\end{table*}

Public UFC is predicted to be the dominant charging method in future smart cities for three reasons. First, their rapid charging times will minimally affect established urban and long-distance driving behavior, thus preserving smooth traffic operation and saving the tremendous cost of renovating traffic control infrastructures. Second, UFC can resolve over-frequent and non-optimal charging behaviors, commonly referred to as ``range anxiety'' \cite{neubauer2014impact}.  This pessimistic attitude towards range highly influences charging and route choice decisions \cite{yang2016modeling}.  In Japan, for example, drivers with access to fast charging tend to come into charge points with a much lower state of charge (SOC) than those without \cite{chang2012financial}. Finally, UFC can greatly improve end-user experience. It is projected that between 3.4\% and 8.5\% of all trips would require fast charging \cite{nicholas2012dc} and drivers believe increased fast charging is essential \cite{azadfar2015investigation}.

The pioneer efforts within UFC deployment are supported by mature technologies. UFC levels have been reached in practice, most notably with regards to light rail transit (LRT), similar to a trolley, which can provide between $200$ and $600kW$, supplying $1.6kWh$ of charge at the bus stop \cite{su2017power, Justino2016Hundreds}. Pilot projects aiming to bring UFC to passenger vehicles have been undertaken worldwide. For example, Porsche installed a prototype $350kW$, $800V$ CS at the their Berlin office in July, 2017. The European Commission's Ultra-E project aims to bring $350kW$ chargers to transportation networks by 2018. Additionally, a rich body of literature from research communities, e.g., power electronics \cite{crosier2012grid, vasiladiotis2015modular} and materials science \cite{kang2009battery}, indicates available hardware and manufacturing techniques are ready for UFC implementation \cite{veneri2013performance, aggeler2010ultra}. 

Despite UFC attractiveness and technological readiness, planning UFC stations in urban areas presents a complex problem, requiring the integration of energy supply, traffic flow, and user utility considerations. The power grid, in particular, will experience unprecedented challenges. Loading stress from AC and DC Levels 1-3 has been shown to lead to significant increments of power losses and voltage deviations \cite{clement2010impact, acha2010effects, maitra2009evaluation}. Further, transferring energy to a vehicle battery can lower the lifetimes of critical grid equipment, e.g., transformers and voltage regulators \cite{Sexauer2013App, yan2012impact, Argade2012Prob}. UFC power density, 200 times that of AC Level 1, could cause detrimental impacts to the power grid, e.g., frequent blackouts, causing great concerns for grid operators and utilities and prohibiting UFC deployment. 

A great amount of research has been devoted to CS planning from both the engineering and social-science communities. The scopes of these studies, however, are limited by existing charging levels (AC and DC Level 1-3) and integration constraints, such as traffic flow and user satisfaction, which does not allow for a direct solution to be drawn for UFC planning. 


For this reason, we present a comprehensive review of the technical components required for the construction of a complete solution for UFC planning. In Section \rom{2} we introduce the charging characteristics of UFC stations and discuss the inherent differences between conventional planning methods and those needed for UFC planning. In Section \rom{3} UFC planning formulations are established with respect to traffic flow and user behavior. To address the needs of the power grid,  Section \rom{4} introduces the analogy between renewable generation and UFC with regards to how their electric characteristics present to the power grid. Based on this result, the methods for determining UFC stations' geographical location (siting) and charging capacities (sizing) are presented with suitable formulations and algorithms from Renewable Energy (RE) planning literature.  Finally, additional considerations and future research directions are presented in Section \rom{5}. The contributions of this paper are twofold:
\begin{itemize}
\item A review of planning methods for PEV charging stations under a wide scope of power levels (AC and DC Level 1-3), constraints (traffic flow and user utility), and research disciplines (engineering and social-science). 
\item An outlook for methods of integrating UFC stations within smart cities under multi-faceted considerations, in particular constraints from the power grid, by drawing analogies with RE planning. 
\end{itemize}

\section{Distinguishing Ultra-Fast Charging Stations}
Currently, the fastest publicly available charging level is ``DC fast charging''. SAE standards \cite{sae2011sae} characterize fast charging (FC) by a DC power level at or less than $40kW$. Thus, a station consisting of 6 sockets could be expected to draw about $240kW$ during a period of full use, providing a full charge within 30 minutes. Within industry implementation, most publicly available DC fast chargers use CHAdeMO standards with an output of $50kW$ or less and Tesla's Superchargers output up to $135kW$. 

Under the ambition of boosting charging speeds to a level comparable to traditional refueling, UFC has been proposed \cite{hoimoja2012toward,vasiladiotis2015modular,aggeler2010ultra}. At present, there is no standardized power rating to define UFC, due to the lack of wide-spread implementation. A working definition is within the range of $350kW$ to $1MW$, per vehicle per charge, connected to 3-phase DC power. In continuing towards the goal of reaching traditional refueling stations, 6 or more sockets per station location are expected, bringing total station power draw to withing $2$ to $6MW$ \cite{crosier2012grid}. 

With such high power requirements, the connection topology and configuration between charging stations and the grid requires adjustment. High power demand may require a high input voltage level, necessitating UFC connections at the transmission/sub-transmission system instead of the traditional medium/low voltage distribution system. Such a configuration is  similar to an industrial load. So far, besides conceptual proposals and a few pilot projects detailed within the Introduction, the feasibility, topology, and potential grid impacts of UFC have yet to be thoroughly studied. 
 
\subsection{Comparisons to Existing Charging Station Planning Methodologies} 
In essence, UFC planning must improve on traditional siting and sizing by considering two central characteristics: (i) stochasticity and (ii) user behavior. Past CS planning research has considered both the power grid and end-users, in the form of transportation networks and/or user behavior.  User behavior and traffic flow manifest in similar ways whether a traditional or UFC  station is being considered, requiring little change. 

The key difference between UFC and traditional siting and sizing appears when the power requirements are considered. For example, when only FC is considered, power requirements for the charging site are around $50kW$, which influences connection locations in the grid \cite{yao2014multi}. UFC, on the other hand, can require up to $1MW$ per vehicle. This difference impacts the parameters of the problem and the system under study, i.e., connections must be made more similarly to an industrial customer. Further, current CS planning research focuses on urban networks \cite{jia2014novel}. This is still of value for UFC planning, but a key contribution of UFC is with respect to long distance trips. Thus, rural areas and highways must be considered, inherently expanding and changing the layout of the networks of interest from low voltage distribution systems (less than $69kV$) to subtransmission ($69-138kV$) and transmission systems (greater than $138kV$). 

The highly impulsive nature of charging loads is not often integrated within planning formulations \cite{Sexauer2013App}, but must be considered with UFC's high power density. Frequently, the stochasticity of charging events is not considered  \cite{wei2017expansion,liu2013optimal}. When determining where to site CS that require much lower power levels, stochasticity is often overlooked, as these charging events do not impact the grid in an exceptional way and the consideration of probabilistic events can further complicate already complex formulations. Many choose to instead generate candidate sites based on service levels \cite{liu2013optimal} or estimated traffic demands \cite{trivedi2015multiobjective} or use static data \cite{liao2016planning}.  For example, general curves indicating parking lot usage can be created to give an estimate of charge usage \cite{trivedi2015multiobjective}. However, in the face of CS that can cause megawatt increases to the grid demand, stochastic behavior must be fully considered in order to determine truly optimal locations and sizes that the grid is capable of servicing, even in times of frequent use. Though there are inherent differences with respect to charging power characteristics, previous methods provide a starting point for UFC planning.  

\subsection{Stochastic and Multistage Considerations}
Much of the aforementioned research ignores the stochastic behavior of loads and/or generation capacity or include probability density functions or other methods that do not fully capture stochastic characteristics. Formulaically,  established two-stage and multi-stage formulations are used to address uncertainty \cite{santos2017new, wang2015stochastic}. Sample average approximation (SAA) can solve two-stage formulations and have been used for power system problems, e.g., generation/transmission expansion and RE usage. SAA will converge to an optimal solution with large sample sizes, which may not be available for CS usage patterns. Multistage formulations typically employ decomposition approaches, e.g., dual decomposition or Bender's decomposition, both of which are easily solved and, in the case of large systems, parallel processing can increase solution speed. \\
\indent Transmission expansion planning commonly employs stochastic or robust formulations to solve system expansion under uncertainties that are attached to extremely high-cost decisions \cite{chen2016robust,qiu2017probabilistic}. Unit retirement is included within these problems, but can be removed in the case of UFC expansion, as no equivalent exists. Regardless, these models show promise for application to charging station placement due to the consideration of uncertainty in size and location of candidate sites. Within two-stage problems, the first stage contains no uncertainty and is made with respect to each possible scenario, e.g., demand, within the second stage. For example, sites and sizes are determined in the first stage and, upon realization of uncertainties, the second stage takes recourse action to address possible infeasibilities. Solution procedures for larger problems typically includes a calculation of a proposed expansion plan with respect to worst case uncertainties and then optimal power flow is calculated based on the proposed plan.

\begin{figure}[!tp]
\centering
\includegraphics[width=3.5in]{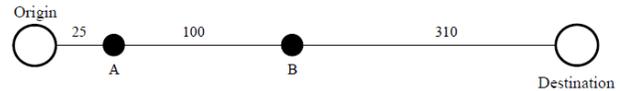}
\caption{A possible OD pair and facility location. For the FCLM, note that for an EV with a range of 350 mi, a facility would be required at node B to ensure a completed trip from origin to destination.  In the case of the FRLM Model, the importance of round trip consideration is apparent, as a vehicle with a range of 700mi would only be able to complete a round trip with a station placed at node B.}
\label{TFEX}
\end{figure}

\section{The Importance of User Driving and Charging Behavior}
Gas station and public CS siting have been studied using algorithmic tools considering two sets of constraints, traffic flow and customer behavior, which are equally applicable to UFC planning with adequate adjustments to customer charging habits.  However, these tools do not address major limitations caused by the power grid. A review of the mathematical core elements for these approaches is presented in this section.
  
\subsection{Traffic Flow Models} 
A key consideration for CS placement, outside of power system characteristics, is the flow of traffic within the region of interest. Most traffic flow models are built from the standard version of location theory and the initial \textit{flow capturing location model (FCLM)} developed in 1990 \cite{hodgson1990flow}, for locating retail facilities. Knowledge of traffic flow between origin and destination (OD) pairs facilities the maximization of flow captured. Refueling at multiple stations cannot occur, a deficiency addressed by the \textit{flow-refueling location model (FRLM)}, a mixed-integer non-linear program (MINLP) \cite{kuby2005flow,capar2012efficient}. If the shortest path from origin to destination and back can be traversed without running out of fuel, flow is satisfied. An example OD pair and two possible facility locations are shown in Fig. \ref{TFEX}. In addition to OD pairs and flow volumes inputs, the FRLM requires pregeneration of OD pair sets and all combinations of possible refueling facilities for each path. 

Facility capacity is not considered in either the FRLM or FCLM, but is included within the \textit{capacitated flow-refueling model (CFRLM)} \cite{upchurch2009model,hosseini2017heuristic}, which requires knowledge of where refueling occurs to ensure capacity is not exceeded. The CFRLM also reduces required constraints by removing an indicator variable for facility combinations and presents as a MINLP, typically solved with branch-and-bound. The pregeneration process needed by the FRLM and FCLM is removed with \textit{arc cover-path cover (ACPC)} methods using candidate sets \cite{capar2013arc}. OD pairs are replaced with paths between two points, with flow considered refueled if all arcs on one path can be traversed. However, the 1:1 correspondence of arcs and constraints results in increasing problem size with increasing arcs.    
 
An extension of the ACPC model, multi-period planning, considers the capacity of stations and demand dynamics with a nonlinear formulation and similar assumptions as previous efforts (shortest path, \textit{a priori} knowledge of demand, etc.)\cite{zhang2017incorp}. Traffic is served if a PEV can go from origin to destination and back without running out of battery, similar to the FRLM. This model is solved via a heuristic method (``Forward Method''), in which periods are solved individually and previous solutions are passed to subsequent periods, yielding a relatively small solution gap. Various other works have presented alterations to these fundamental formulations and most use heuristic approaches to accomodate problem structure and size \cite{ghamami2016general, arslan2016benders, sathaye2013approach, wang2009locating}.

Some have also attempted to bridge the gap between transportation planning and the power system characteristics of charging stations \cite{wang2013traffic, wei2017expansion}. Pairing power system and traffic flow requires alterations to objectives and constraints in the face of two conflicting objectives \cite{zhang2016pev, yao2014multi}. These methods do not consider patterns of CS usage, which are highly stochastic, and assume known traffic patterns.

\subsection{User Behavior and Needs} 
Many of the user behavior studies that have been conducted use on-board trip logging and surveys \cite{azadfar2015investigation,tal2014charging,vilimek2012mini,  franke2013understanding}. Though these studies shed light on consumer behavior, they fall victim to the relatively small amount of EVs on the road. CS themselves can also be monitored, indicating how long and where drivers charge \cite{morrissey2016future}. These studies can give information about accurate usage assumptions within case studies. For example, assuming a driver only charges from empty to full charge may not be an accurate. Studies can also divide CS not only based on geographic location, but by what is near them within the city, e.g., on-street, retail locations, gas stations, etc.  Using these general trends, weighting factors can be determined for candidate locations based on their surroundings \cite{wagner2013optimal}.

Statistical approaches can be used to determine the best estimators of user behavior. Some studies create only probability distribution functions of the data in order to create representative profiles \cite{quiros2015statistical}. An interesting approach to modeling PEV user behavior uses of past driver data for traditional ICE vehicles, i.e., it is assumed that drivers will conduct the same trips \cite{kelly2012time}. Logit-models, frequently used for models of repeated choice, are used to predict behavior based on generalized profiles \cite{luo2015consumer} and can determine behavioral output, elasticity of choice, and valuation of attributes \cite{jabeen2013electric}. Regardless, a glaring deficiency exists with regards to data that properly characterizes charging behavior, which is critical to determining high value locations.  

\section{Renewable Siting and Sizing for UF Charging Station Planning}
When the grid interacts with UFC technology, it sees a new type of electric load, which possess three distinct features, (i) an extremely high power, impulsive demand, (ii) inherently difficult to predict charging activity, and (iii) a power-electronics-dominated high-bandwidth control framework. Interestingly, these features are very similar to aspects studied within the field of RE integration with respect to stochasticity and high power requirements. For example, a typical moderately sized PV installation can reach an average output of $5-9MW$, comparable to a UFC station \cite{kumar2015performance}. Further, RE output highly depends on environmental factors, which are stochastic in nature. PV installations, for example, can see production drop to 10-25\% of rated capacity on cloudy days. Charging behavior, dominated by human choice, is also highly stochastic. The high power stochastic qualities of RE and charging behavior can both harm grid performance. CS impacts have been detailed above and, similarly, RE can cause transient and voltage instability and power quality issues \cite{kou2017PQ}. Inspired by the similarities between these high powered, intermittent technologies, and given the more plentiful research focused on RE planning, it is helpful to begin looking at the stochastic UFC load as a type of reverse RE source. 

Siting and sizing is a key consideration for both RE and UFC planning. Approaches for the integration of RE within the grid can generally be classified into three different types, mathematical, heuristic, and numerical, all of which are summarized in Table \ref{methods}. These three approaches can be modified and applied to UFC in different ways, with power system characteristics serving as hard constraints. These methods yield insight into the impact of various sites and sizes of RE on power grid performance, a facet of analysis also of interest for those planning UFC expansion.

\begin{table*}[!t]
\begin{center}
\caption{Renewables Siting and Sizing Methodologies} \label{methods} 
    \begin{tabular}{ | c | c | c | c |}
    \hline
    Approach & \thead{Solution Methods} & \thead{Benefits} & \thead{Concerns} \\ \hline 
    Mathematical & \thead{Mixed Integer Linear/Nonlinear Program \cite{arriaga2016long,singh2015allocation,al2013optimal,yang2015milp} \\ Quadratic Program \cite{sfikas2015simultaneous} \\ Optimal Power Flow \cite{mahmoud2016optimal}} & \thead{Existing Solvers \\ Established Formulations} & \thead{Computationally Expensive \\ Data and Situation Specific} \\ \hline
   Heuristic & \thead{Genetic Algorithms \cite{sheng2015optimal, mena2014risk} \\ (Open - Space) Particle Swarm \cite{abdul2012optimum} \\ Harmony Search \cite{nekooei2013improved} \\ Artificial Bee Colony \cite{abu2011optimal,sebaa2014optimal}} &\thead{Speed of Solution \\ Ease of Use} & \thead{Global Optimal Not Guaranteed \\ Susceptible to Parameter Tuning } \\ \hline
    Numerical & \thead{AC Power Flow \cite{mahmoud2016optimal} \\ Exhaustive Power Flow \cite{elsaiah2014analytical} \\ Modified Power Flow \cite{acharya2006analytical,griffin2000placement}} & \thead{Relatively Simple Problem Structure \\ Ease of Formulation} & \thead{Only Power System Characteristics \\ Requires Input Demand Profile \\ No Separation of Site and Size} \\ \hline
    \end{tabular}
\end{center}
\end{table*}

\subsection{Mathematical Approaches}
Mathematical approaches are those methods based on the use of established mathematical programs to solve the issue of RE integration with respect to established problems, e.g., those with known or specific data. Commonly implemented formulations are mixed integer linear programs (MILP) \cite{yang2015milp}, mixed integer nonlinear programs (MINLP) \cite{arriaga2016long,singh2015allocation,al2013optimal} , quadratic programs (QP) \cite{sfikas2015simultaneous}, and AC optimal power flow based procedures \cite{mahmoud2016optimal}. Formulations are built from the ground up using system characteristics to build operational and problem-specific constraints.  Due to the nature of power systems, many of these formulations are easily modified to meet the requirements of locales different than their original formulation. When faced with known data, e.g., weather patterns, heirarchical decisions can be made to limit candidate sites and create a priority list \cite{rumbayan2012prioritization}, enabling quicker and easier solution of the overall MINLP \cite{arriaga2016long}. Additional input data typically includes geographical locations of aspects of the systems (e.g., buildings or available sites) \cite{arriaga2016long}, load profiles \cite{singh2015allocation}, and electricity pricing \cite{arriaga2016long}. 

Formulations of this type are easily handled with solvers in existing packages, e.g., MATLAB, GAMS, CPLEX, and their convergence to optimal solutions can be proven with optimization theory, e.g., duality or convexity. Non-linear formulations can also be convexified, using mathematical modifications \cite{wei2017expansion}. The critical issue with these implementations comes in the form of overly complex constraints, which can be represented in two forms: (i) highly nonlinear constraints induced by renewable output \cite{al2013optimal} (similar to charging behavior) and power grid dynamics, e.g., thermal loading and generation unit frequency \cite{sfikas2015simultaneous,Vallem2005siting}, and (ii) high numbers of constraints introduced by the consideration of large-scale problems, also likely to occur when considering traffic flow within UFC planning \cite{yao2014multi, wei2017expansion}. These complexities affect the convergence properties of the planning problem, but can be addressed, e.g., non-linearities can be linearized or convexified \cite{wei2017expansion, fiorini2017sizing}. Many exclude highly non-linear aspects, e.g., harmonics, thermal loading, and admittances, completely, which may limit the application of solutions to real world problems. These potential complications indicate the importance of creating valuable formulations and utilizing mathematical strategies to overcome burdensome constraints. Sparse structures, for example, have been used to aide in the solution of complex formulations \cite{liu2013optimal}.  In order to address the computational burden of these methods, many have turned to heuristic approaches. 

\subsection{Heuristic Approaches}
Heuristic approaches to RE planning are commonly based on population-based optimization methods, i.e., advanced artificial intelligence algorithms. These approaches minimize losses within the system or otherwise ensure satisfactory performance. Many include objectives concernes with costs of installation/maintenance, active and reactive power, and unserved demand \cite{abdul2012optimum, mardaneh2004siting}.
Genetic algorithms (GA) are widely implemented \cite{sheng2015optimal, mena2014risk}, but can yield results that are of lower quality than those obtained by particle swarm optimization (PSO) techniques \cite{singh2009multiobjective}. PSO can be further improved by extending the search space and assuring reliability of the result, achieved via open-space particle swarm optimization \cite{abdul2012optimum}. Benefits of PSO include ease of coding and overall simplicity, ease of use, high convergence, and a relatively small storage space \cite{ameli2014multiobjective,aman2014new}. Convergence to a global optimum, however, is far from guaranteed and is a pitfall of all heuristic techniques. Pre-processing can be applied to limit the search space \cite{mardaneh2004siting}. 

Artificial Bee Colony (ABC) methods is meta-heuristic method that imitates the foraging behavior of a honeybee swarm to obtain an optimal solution to a constrained or unconstrained problem \cite{abu2011optimal,sebaa2014optimal}.  The results of this method have been shown to perform as well or better than other heuristic approaches \cite{abu2011optimal}, and some implementations require the tuning of only two parameters to obtain the most efficient solution. Further, this heuristic is capable of handling MINLPs, commonly found in CS planning problems. Harmony search \cite{nekooei2013improved} and big bang crunch \cite{sedighizadeh2014application} can also be utilized for RE planning problems. 

Irrespective of the heuristic method being employed, many include multi-objective functions and constraints considering line loss \cite{abu2011optimal}, voltage deviation \cite{abu2011optimal, sebaa2014optimal, nekooei2013improved}, load profiles \cite{singh2009multiobjective}, voltage stability margins \cite{sedighizadeh2014application} , and other power system characteristics, all of which apply within UFC planning. Heuristic methods address complex, multi-objective formulations, which mathematical formulations may be unable to do.  Convergence, however, is highly dependent upon the tuning of critical parameters and is not guaranteed \cite{abu2011optimal}. Heuristic methods often yield solution gaps larger than those of mathematical methods \cite{foster2014comparison}.  

\subsection{Numerical Approaches}
Numerical approaches seek to minimize system loss or improve the voltage profile of the system and have been applied to both radial and meshed networks \cite{wang2004analytical}.  Many of these approaches are based on the well-developed exact loss formula \cite{elgerd1982electric}, which relates voltage magnitude and voltage angle at a bus with the active and reactive power injections (in a highly nonlinear way). These approaches are less complicated than heuristic methods, but require the solving of power flow for exhaustive amounts of case studies. AC power flow yields higher precision, but at a higher computational burden, leading to common implementation of DC power flow instead. There have been improvements to exhaustive search methods and attempts to capture system characteristics as in AC power flow, yielding almost exactly the same results as those obtained via exhaustive power flow methods. Sensitivity analyses, which linearize a nonlinear equation around the initial operation point, are called upon to select appropriate locations in a type of ``pre-processing'', reducing the search space. Pre-generation methods are commonly employed \cite{acharya2006analytical,griffin2000placement} and one such work developed a method for creating a priority list based on power loss sensitivity factors to place and size DG systems for loss reduction \cite{elsaiah2014analytical}. This approach is easily applied to UFC station installation, e.g., representative flows are calculated to determine those locations corresponding to the lowest grid impact and sizing can be determined using the priority list. Additionally, this method can be used on unbalanced networks  by altering the admittance matrix. A relatively simplistic solution structure and focus on loss reduction or voltage profile management makes these procedures attractive to UFC planners \cite{hung2013multiple,hung2010analytical}. The pre-processing stage could be used to select sites corresponding to existing gas stations (assuming similar behavior for ICE and PEVs), addressing user behavior within a power system heavy analysis.

A critical downfall of numerical approaches is their dependence on the demand profile input to the solver and that solutions cannot accurately represent the separate variables of site and size for either renewable generation or CS. Further, many of the studies employing these approaches focus on siting and sizing only one renewable unit at a time. Some have attempted to conduct simultaneous planning for multiple units, with high fidelity results \cite{santos2017new}.

\section{Future Work \& Discussion}
Plug-in electric vehicle market shares have continued to grow in the past seven years, a trend projected to continue with countries pushing towards adoption goals. Various studies have shown that availability of charging speed and available infrastructure play a critical role in PEV adoption rates. Though past studies have focused on the CS placement problem, optimization problems focus on existing charging levels. The pioneer vision for PEV charging, UFC infrastructure, supplies from $350kW$ to $1 MW$ per charge and guides the concept of PEV charging towards speeds comparable to traditional gas refueling. the planning of these UFC stations must consider not only power system requirements, but also user behavior and needs. To address the high power, impulsive characteristics of these new charging loads, RE planning strategies, which consider similar stochastic and power grid related characteristics, can be adopted and combined with existing traffic flow and user behavior models. Regardless, the deficiency on user behavior information must be addressed.

Future UFC work can consider on-site storage and energy management through the lens of industrial customers, comparable in size to UFC stations. Industrial analogies are most beneficial with respect to extensions of UFC, specifically in terms of operation. Much research focusing on industrial loads addresses energy management and, recently, the use of on-site generation to offset the cost of electrical consumption for industrial customers and increase the reliability of their power supply \cite{pipattanasomporn2005implications}. Due to high power requirements, on-site generation has been a topic of research for smaller scale CS \cite{guo2014stochastic}. The participation of industrial customers within curtailment and demand side management programs can also be investigated with an eye towards integration within CS infrastructure, due to similar operational load sizes. The interuptability of CS loads must be studied. The most promising solution for managing CS loads comes in the form of demand side management, by sending favorable ``fuel'' prices at times that benefit grid operation \cite{shojaabadi2016optimal}. Energy management strategies for both industrial customers and smaller scale charging stations have been widely studied \cite{guo2014stochastic, stluka2011energy,ruangpattana2011optimization,pipattanasomporn2005implications}. Due to the similar sizes of these facilities and the proposed CS facilities, energy management schemes will be incredibly valuable. Though these strategies are out of the scope of this paper, they are of key interest for UFC operation.

\bibliographystyle{IEEEtran}

\bibliography{References}

%
%
%

\end{document}